\documentclass{amsart}
\usepackage[all]{xy}
\usepackage{epsf}
\evensidemargin=\oddsidemargin
\newtheorem{thm}{Theorem}[section]
\newtheorem{prop}[thm]{Proposition}

\newtheorem{defn}[thm]{Definition}
\newtheorem{ex}[thm]{Example}

\theoremstyle{definition}	

\numberwithin{equation}{dummy}
\newcommand{\br}{\mathbb R}

\newcommand{\bz}{\mathbb Z}

\newcommand{\mk}{\vskip .1in}

\theoremstyle{plain}

\begin{document}

\title{A survey of Wall's finiteness obstruction}
\author{Steve Ferry}
 \address{S.F.: Dept. of Mathematics\\
 Rutgers University\\ 
 New Brunswick, NJ 08903, USA}
\email{sferry@math.rutgers.edu}

 \author{Andrew Ranicki}
 \address{A.R.: Dept. of Mathematics and Statistics\\ 
 University of Edinburgh\\ 
 Edinburgh EH9 3JZ\\ 
 Scotland, UK}
 \email{aar@maths.ed.ac.uk}


\maketitle

\section*{Contents}

\begin{itemize}
\item[~] Introduction
\item[1.] Finite domination
\item[2.] The projective class group $K_0$
\item[3.] The finiteness obstruction
\item[4.] The topological space-form problem
\item[5.] The Siebenmann end obstruction
\item[6.] Connections with Whitehead torsion
\item[7.] The splitting obstruction
\item[8.] The triangulation of manifolds
\item[~] References
\end{itemize}

\section*{Introduction}

Wall's finiteness obstruction is an algebraic $K$-theory invariant
which decides if a finitely dominated space is homotopy equivalent to a
finite $CW$ complex.  The invariant was originally formulated in the
context of surgery on $CW$ complexes, generalizing Swan's application
of algebraic $K$-theory to the study of free actions of finite groups
on spheres.  In the context of surgery on manifolds, the invariant
first arose as the Siebenmann obstruction to closing a tame end of a
non-compact manifold.  The object of this survey is to describe the
Wall finiteness obstruction and some of its many applications to the
surgery classification of manifolds.  The book of Varadarajan \cite{Va}
and the survey of Mislin \cite{Mis} deal with the finiteness
obstruction from a more homotopy theoretic point of view.

\section{Finite domination}

A space is finitely dominated if it is a homotopy retract of a finite
complex. More formally:

\begin{defn}{\normalfont
	A topological space $X$ is \textit{finitely dominated} 
	if there exists a finite $CW$ complex $ K $ with maps 
	$ d:K \to X $, $s:X \to K $ and a homotopy 
	$$ d \circ s\simeq 	\text{id}_X:X \to X ~.\eqno{\qed}$$}
	\end{defn}

\begin{ex} {\normalfont 
(i) A compact $ANR$ $X$ is finitely dominated (Borsuk \cite{Bo1}). 
In fact, a finite dimensional $ANR$ $X$ can be embedded 
in $\br^N$ ($N$ large), and $X$ is a retract of an
open neighbourhood $U \subset \br^N$ -- there exist a retraction $r:U \to X$
and a compact polyhedron $K \subset U$ such that $X \subset K$,
so that the restriction $d=r\vert:K \to X$ and the inclusion $s:X \to K$ 
are such that $d\circ s=\text{id}_X:X \to X$ .\\
(ii) A compact topological manifold is a compact $ANR$, and hence finitely
dominated.}\hfill\qed
\end{ex}

The problem of deciding if a compact $ANR$ is homotopy equivalent to
a finite $CW$ complex was first formulated by Borsuk \cite{Bo2}.
(The problem was solved affirmatively for manifolds by Kirby and
Siebenmann in 1969, and in general by West in 1974  -- see section
8 below.) The problem of deciding if a finitely dominated space is
homotopy equivalent to a finite $CW$ complex was first formulated
by J.H.C.Whitehead.
Milnor \cite{Mil1} remarked: ``It would be interesting to ask if every space
which is dominated by a finite complex actually has the homotopy type of a
finite complex. This is true in the simply connected case, but seems like a 
difficult problem in general.''

\medskip

Here is a useful criterion for recognizing finite domination:

\begin{prop}
A $CW$ complex $X$ is finitely dominated if and only if 
there is a homotopy $ h_{t}:X \to X $ such that $ h_0=id $ and 
$ h_1(X) $ has compact closure.
\end{prop}
\begin{proof} If $ d:K \to 
X $ is a finite domination with right inverse $ s $, let $ h_{t} 
$ be a homotopy from the identity to $ d \circ s $.  Since $ 
h_1(X) \subset d(K) $, the closure of $ h_1(X) $ is compact in 
$X$.  Conversely, if the closure of $ h_1(X) $ is compact in 
$X$, let $ K $ be a finite subcomplex of $X$ containing $ 
h_1(X) $.  Setting $d$ equal to the inclusion $ K \to X $ and 
$ s $ equal to $ h_1:X \to K $ shows that $X$ is finitely 
dominated.
\end{proof}

It is possible to relate finitely dominated spaces, finitely
dominated $CW$ complexes and spaces of the homotopy type of $CW$ complexes,
as follows.

\begin{prop}
{\rm (i)} A finitely dominated topological space $X$ is homotopy equivalent to a
countable $CW$ complex.\\
{\rm (ii)} If $X$ is homotopy dominated by a finite $ k $-dimensional $CW$ 
	complex, then $X$ is homotopy equivalent to a countable $ (k+1) 
	$-dimensional $CW$ complex.
\end{prop}

\begin{proof} The key result is the trick of Mather \cite{Ma}, which shows 
that if $d:K \to X$,
$s:X \to K$ are maps such that $d\circ s \simeq \text{id}_X:X \to X$ then $X$ is 
homotopy equivalent to the mapping telescope of $s \circ d :K \to K$.
This requires the calculus of mapping cylinders, which we now recall.\\
By definition, the mapping cylinder of a map $f:K \to L$ is the identification space
$$M(f)~=~(K \times [0,1] \cup L)/((x,1) \sim f(x))~.$$
We shall use three general facts about mapping cylinders:
	\begin{itemize}
		\item  If $ f:K \to L $ 
		and $ g:L \to M $ are maps and $ k:K \to M $ is homotopic to $ 
		g \circ f $, the mapping cylinder $ M(k) $ is homotopy 
		equivalent rel $ K \cup M $ to the concatenation of the mapping 
		cylinders $ M(f) $ and $ M(g) $ rel $ K \cup M $.
		\item  If  $ f,g:K \to L $ with $ f \sim g $, then the mapping 
		cylinder of $ f $ is homotopy equivalent to the mapping cylinder 
		of $ g $ rel $ K \cup L $.
		\item  Every mapping cylinder is homotopy equivalent to its base rel 
		the base.
	\end{itemize}

\noindent The mapping telescope of a map $\alpha:K \to K$ is the countable union
$$\bigcup_{i=0}^{\infty}M(\alpha)~=~
\bigcup_{i=0}^{\infty} K \times [i,i+1]/\{(x,i) \sim (\alpha(x),i+1)\}~.$$
For any maps $d:K \to X$, $s:X \to K$ we have
$$\bigcup_{i=0}^{\infty}M(d \circ s) ~=~ X \times I \cup \bigcup_{i=0}^{\infty}M(s \circ d)$$ 
with $\bigcup_{i=0}^{\infty}M(s \circ d)$ a deformation retract, so that
$$\bigcup_{i=0}^{\infty}M(d \circ s) ~\simeq~ \bigcup_{i=0}^{\infty}M(s \circ d)~.$$
To see why this holds, note that 
	\( \bigcup_{i=0}^{\infty}M(d \circ s) \) is homotopy 
	equivalent to an infinite concatenation of alternating \( M(d) \)'s and \( M(s) 
	\)'s which can also be thought of as an infinite concatenation of \( 
	M(s) \)'s and \( M(d) \)'s. Essentially, we're reassociating an 
	infinite product. Here is a picture of this part of the 
	construction.
	

\mk
\centerline{\epsfysize=60mm\epsfbox{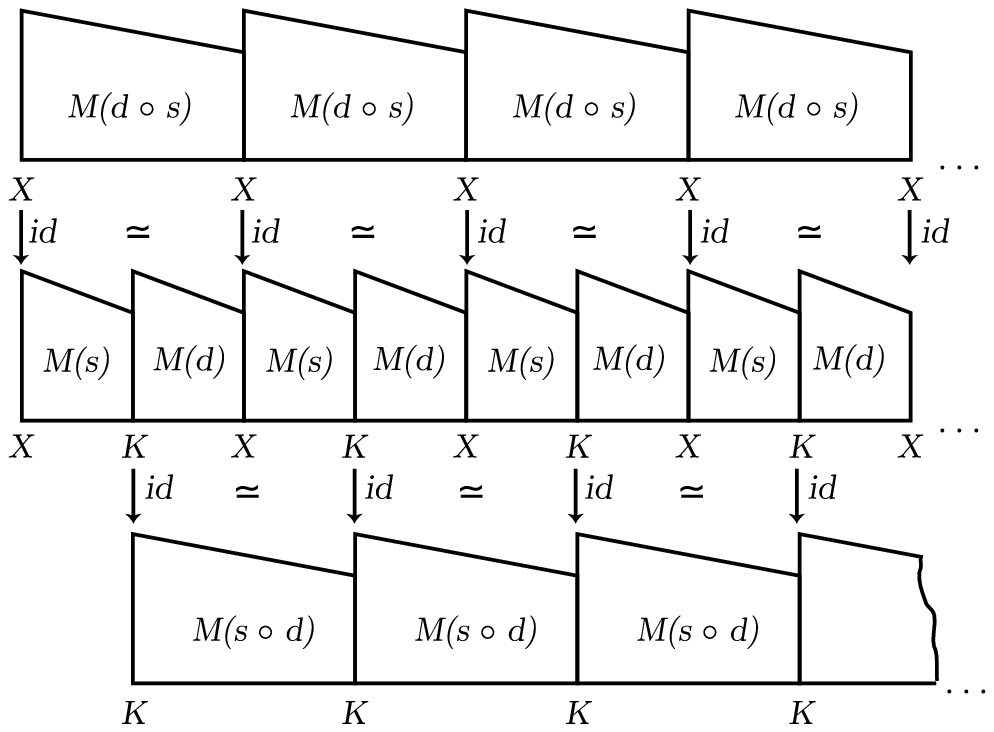}}
\mk


\noindent (i) If $d:K \to X $,  $s:X \to K$ are such that $d \circ s \simeq \text{id}_X:X \to X$
there is defined a homotopy idempotent of a finite $CW$ complex
$$\alpha~=~s \circ d:K \to K~,$$
with $\alpha \circ \alpha \simeq \alpha:K \to K$.
We have	homotopy equivalences
	\[ X \simeq X \times [0,\infty) \simeq \bigcup_{i=0}^{\infty}M(\text{id}_{X}) \simeq 
	\bigcup_{i=0}^{\infty}M(d \circ s) \simeq \bigcup_{i=0}^{\infty}
	M(s \circ d) = \bigcup_{i=0}^{\infty}M(\alpha). \]
The mapping telescope $\bigcup_{i=0}^{\infty}M(\alpha)$ is a countable 
$CW$ complex.\\
(ii) As for (i), but with $K$ $k$-dimensional.
\end{proof}

This proposition is comforting because it shows that the finiteness 
problem for arbitrary topological spaces reduces to the finiteness 
problem for $CW$ complexes.  One useful consequence of this is that we 
can use the usual machinery of algebraic topology, including the 
Hurewicz and Whitehead theorems, to detect homotopy equivalences.
\mk

\begin{prop} {\rm (Mather \cite {Ma})} \label{Mather}
A topological space $X$ is finitely dominated if and
only if $X \times S^1$ is homotopy equivalent to a finite $CW$ complex.
\end{prop}
\begin{proof}
The mapping torus of a map $\alpha:K \to K$ is defined (as usual) by
$$T(\alpha)~=~(K \times [0,1])/\{(x,0) \sim (\alpha(x),1)\}~.$$
For any maps $ d:K \to X $, $s:X \to K $ there is defined a homotopy
equivalence
$$T(d\circ s:X \to X) \to T(s\circ d :K \to K )~;~(x,t) \mapsto (s(x),t)~.$$
If $d \circ s \simeq \text{id}_X:X \to X$ and $K$ is a finite $CW$ complex
we thus have homotopy equivalences
$$X \times S^1 \simeq T(\text{id}_X) \simeq T(s \circ d)$$
with $T(s \circ d)$ a finite $CW$ complex.\\
Conversely, if $X \times S^1$ is homotopy equivalent to a finite $CW$
complex $K$ then the maps
$$d : K \simeq X \times S^1 \stackrel{proj.}{\longrightarrow} X~~,~~
s : X \stackrel{incl.}{\longrightarrow} X \times S^1 \simeq K$$
are such that $d \circ s \simeq \text{id}_X$, and $X$ is dominated by $K$.
\end{proof}

\section{The projective class group $K_0$}

Let $\Lambda$ be a ring (associative, with 1).

\begin{defn}{\normalfont
	A $\Lambda$-module $P$ is \textit{f.\,g. projective} 
	if it is a direct summand of a f.\,g. (= finitely generated)
	free $\Lambda$-module $\Lambda^n$, with $P \oplus Q = \Lambda^n$
    for some direct complement $Q$.}\hfill\qed
	\end{defn}

    A $\Lambda$-module $P$ is f.\,g. projective if and only if
	$P$ is isomorphic to $\text{im}(p)$ for some projection
	$p=p^2:\Lambda^n \to \Lambda^n$.

\hyphenation{Grothen-dieck}
\begin{defn}{\normalfont
	(i) The \textit{projective class group} 
	$ K_0(\Lambda) $ is the Grothendieck group of stable isomorphism
	classes of f.\,g. 	projective $\Lambda $-modules.\\
	(ii) The  \textit{reduced projective 
	class group} $\widetilde K_0(\Lambda) $ is the quotient of $ 
	K_0(\Lambda) $ by the subgroup generated by formal differences 
    $[\Lambda^m]-[\Lambda^n]$ of 
	f.\,g. free modules.}\hfill\qed
\end{defn}

  Thus an element of 
	$\widetilde K_0(\Lambda) $ is an equivalence class $[P]$ of
	f.\,g. projective $\Lambda $-modules, with
    $ [P_1]=[P_2] $ if 
	and only if there are f.\,g. free $\Lambda $-modules $ 
	F_1 $ and $ F_2 $ so that $ P_1 \oplus F_1 $ is isomorphic to 
	$ P_2 \oplus F_2 $.  In particular, $ [P] $ is trivial if and only 
	if $ P $ is \textit{stably free}, that is, if there is a f.\,g.
	free module $ F $ so that $ P \oplus F $ is free.

\mk
Chapter 1 of Rosenberg \cite{Ro} is a general introduction to the 
projective class groups $K_0(\Lambda)$, $\widetilde{K}_0(\Lambda)$
and their applications, including the Wall finiteness obstruction.
\mk

\begin{ex}{\normalfont
There are many groups $\pi $ for which 
$$\widetilde K_0(\bz[ \pi])~=~0~,$$
including virtually polycyclic groups, a class which includes free and
free abelian groups.}\hfill\qed
\end{ex}  

At present, no example is known of a torsion-free infinite group $\pi$ 
with
$\widetilde{K}_0(\bz[ \pi])\neq 0$.  Indeed, Hsiang has conjectured that
$\widetilde K_0(\bz[ \pi])=0 $ for any torsion-free group $\pi$.  
(See Farrell and Jones \cite{FJ}, pp. 9--11).
On the other hand:

\begin{ex} {\normalfont 
(i) There are many finite groups $\pi $ for which 
$$\widetilde K_0(\bz[\pi])~\neq~0~,$$
including the cyclic group $\bz_{23} $.\\
(ii) The reduced projective class group of the quaternion group 
$$Q(8)~=~\{\pm 1,\pm i, \pm j,\pm k\}$$ 
is the cyclic group with 2 elements
$$\widetilde{K}_0(\bz[Q(8)])~=~\bz_2~,$$
generated by the f.\,g. projective $\bz[Q(8)]$-module 
$$P~=~\text{im}\bigg(\begin{pmatrix}
1-8N & 21N \\[1ex] -3N & 8N
\end{pmatrix}~:~\bz[Q(8)]\oplus\bz[Q(8)] \to
\bz[Q(8)]\oplus \bz[Q(8)]\bigg)$$
with $N=\sum\limits_{g \in Q(8)}g$.}\hfill\qed
\end{ex}

We refer to Oliver \cite{O} for a survey of the computations of 
$\widetilde{K}_0(\bz[\pi])$ for finite groups $\pi$.

\section{The finiteness obstruction}

Here is the statement of Wall's theorem.

\begin{thm} {\rm (\cite{Wa1},\cite{Wa2})}
    {\rm (i)} A finitely dominated space $X$ has a finiteness obstruction 
	$$[X] \in \widetilde K_0(\bz[\pi_1(X)]) $$ 
	such that $[X]=0 $ if and only if $X$ is 
	homotopy equivalent to a finite $CW$ complex.\\
	{\rm (ii)} If $\pi $ is a finitely presented group then every element 
	$\sigma \in \widetilde K_0(\bz[ \pi ]) $ is the finiteness obstruction
	of a finitely dominated $CW$ complex $X$ with $[X]=\sigma $, $\pi_1(X)=\pi$.\\
	{\rm (iii)} A $CW$ complex $X$ is finitely dominated if and only if $\pi_1(X)$
is finitely presented and the cellular $\bz[\pi_1(X)]$-module chain 
complex $C_*(\widetilde X)$ of the universal cover $\widetilde X$
is chain homotopy equivalent to a finite chain complex ${\mathcal P}$
of f.\,g. projective  $\bz[\pi_1(X)]$-modules.
\end{thm}

\noindent{\it Outline of proof} 
(i) Here is an extremely condensed sketch of Wall's argument from 
\cite{Wa1}. If $ d:K \to X $ is a finite domination with $X$ a 
$CW$ complex, we can assume that $d$ is an inclusion by replacing $X$, 
if necessary, by the mapping cylinder of $d$.  For each $ 
\ell \ge 2 $, we then have a split short exact sequence of abelian 
groups 
$$0 \to \pi_{\ell+1}(X,K) \to \pi_{\ell}(K) \to \pi_{\ell}(X) \to 0~.$$ 
Wall gives a special argument to show that $d$ can be 
taken to induce an isomorphism on $\pi_1 $ and then shows that $ 
\pi_{\ell+1}(X,K) $ is f.\,g. as a module over $\bz[ 
\pi_1(X) ]$, provided that $\pi_{q}(X,K)=0 $ for $ q \le \ell,\ 
\ell \ge 2 $.  
This allows him to attach $\ell+1 $-cells to form a complex $\overline 
K \supset K $ and a map $\overline d: \overline K \to X $ extending $d$ 
so that $\overline d $ induces isomorphisms on homotopy groups through 
dimension $\ell $.  Since $\overline d $ is a domination with the 
same right inverse $ s $, this process can be repeated.  In the case $ 
\ell \ge \dim(K) $, Wall shows that $\pi_{\ell+1}(X,K) $ is a 
f.\,g. \textit{projective} module over $\bz[ \pi_1(X)]$.
If $\pi_{\ell+1}(X,K)$ is free (or even stably free) we can attach $\ell+1 
$-cells to kill $\pi_{\ell+1}(X,K) $ without creating new problems 
in higher dimensions.  The result is that $\overline d $ is a homotopy 
equivalence from $\overline K $ to $X$.  If this module is not stably 
free, we are stuck and the finiteness obstruction 
is defined to be 
$$[X]~=~(-1)^{\ell+1}[\pi_{\ell+1}(X,K)]\in
\widetilde K_0(\bz[ \pi_1(X)])~.$$
(ii) Given a finite $CW$ complex $ K $ and a nontrivial $\sigma \in 
\widetilde K_0(\bz[ \pi_1(K)]) $, here is one way to construct a $CW$ 
complex with finiteness obstruction $\pm \sigma $: let $\sigma $ 
be represented by a f.\,g. projective module $ P $ and 
let $F= P \oplus Q $ be free of rank $n$. Let $ A $ be the 
matrix of the projection $p: F \to P \to F $ with respect to a standard 
basis for $ F $. Now let 
\[ L=K \vee \bigvee_{i=1}^nS_{i}^{\ell}. \]
\mk\noindent
There is a split short exact sequence
\[  
\xymatrix{
0 \ar[r] & \pi_{\ell}(K) \ar[r]_{i_{*}} & \pi_{\ell}(L) \ar@<-1ex>[l]_{r_{*}} 
\ar[r] & \pi_{\ell}(L,K) \ar[r] & 0,
}
\]
where $ r:L \to K $ is the retraction which sends the $S^{\ell} 
$'s to the basepoint. Since
\[ 
\pi_{\ell}(L,K) \cong \pi_{\ell}(\tilde L,\tilde K) \cong H_{\ell}(\tilde L,
\tilde K) \cong F,
\]
we can define $\alpha:L \to L $ so that $\alpha|K = id $ and so
that $\alpha_{*}:\pi_{\ell}(L) \to \pi_{\ell}(L) $ has the 
matrix 
\[ \left(
\begin{array}{cc}
	id & 0  \\
	0 & A
\end{array} \right)
 \]
with respect to the direct sum decomposition $\pi_{\ell}(L) \cong 
\pi_{\ell}(K) \oplus F $.  Since $ A^2=A $, it is easy to check 
that $\alpha $ is \textit{homotopy idempotent}, i.e. that $ 
\alpha \circ \alpha \sim \alpha $ rel $ K $. 
\mk
Let $X$ be the infinite direct mapping telescope of $\alpha $ pictured 
below.

\mk
$$\epsfxsize=\hsize \epsfbox{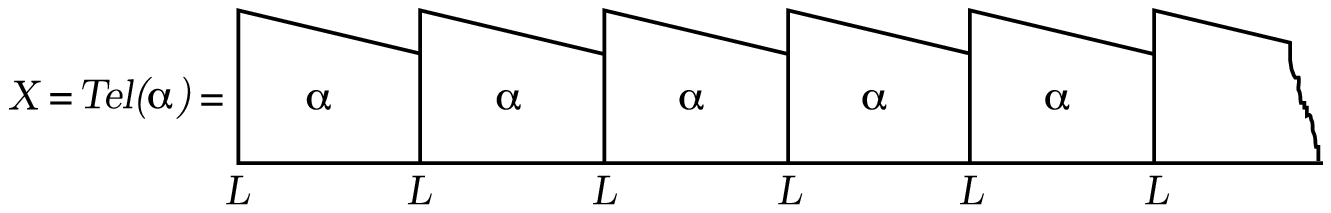}$$
\mk
\mk\noindent
Let $ d:L \to X $ be the inclusion of $ L $ into the top level of 
the leftmost mapping cylinder of $X$ and define $ s':X \to L $ 
by setting $ s' $ equal to $\alpha $ on each copy of $ L $ and 
using the homotopies $\alpha \circ \alpha \sim \alpha $ to extend 
over the rest of $X$. One sees easily that $ d \circ s' $ 
induces the identity on the homotopy groups of $X$ and is 
therefore a homotopy equivalence. If $\phi $ is a homotopy inverse 
for $ d \circ s' $, we have $ d \circ  s \sim id $, where $ s = 
s'\circ \phi $. This means the $d$ is a finite domination with 
right inverse $ s $. It turns out that $[X]=(-1)^{\ell+1}[P] $.
\mk
\noindent
(iii) If $X$ is dominated by a finite $CW$ complex $K$ then $\pi_1(X)$ 
is a retract of the finitely presented group $\pi_1(K)$, and is thus also
finitely presented. The cellular chain complex $C_*(\widetilde{X})$ is a
chain homotopy direct summand of the finite f.g. free $\bz[\pi_1(X)]$-module
chain complex $\bz[\pi_1(X)]\otimes_{\bz[\pi_1(K)]}C_*(\widetilde{K})$, with
$\widetilde{K}$ the universal cover of $K$. It follows from the algebraic
theory of Ranicki \cite{Ra1} (or by the original geometric argument of
Wall \cite{Wa2}) that $C_*(\widetilde{X})$ is  chain equivalent
to a finite f.g. projective $\bz[\pi_1(X)]$-module chain complex
${\mathcal P}$.
\mk\noindent 

Conversely, if $\pi_1(X)$ is finitely presented and $C_*(\widetilde{X})$
is  chain equivalent to a finite f.g. projective $\bz[\pi_1(X)]$-module 
chain complex ${\mathcal P}$ the cellular $\bz[\pi_1(X)][z,z^{-1}]$-module
chain complex of the universal cover $\widetilde{X \times S^1}=
\widetilde{X}\times \br$ of $X \times S^1$ 
$$C_*(\widetilde{X \times S^1})~=~C_*(\widetilde{X})\otimes_{\bz}C_*(\br)$$
is chain equivalent to a finite f.g. free $\bz[\pi_1(X)][z,z^{-1}]$-module
chain complex, so that $X \times S^1$ is homotopy equivalent to a finite
$CW$ complex (by the proof of (i), using $[X \times S^1]=0$) and $X$ is 
finitely dominated by 1.5.\hfill \qed
\mk

In particular, if $\pi$ is a finitely presented group such that 
$\widetilde{K}_0(\bz[\pi])\neq 0$ then there exists a finitely dominated
$CW$ complex $X$ with $\pi_1(X)=\pi$ and such that $X$ is not homotopy
equivalent to a finite $CW$ complex. See Ferry \cite{Fe1} for the
construction of finitely dominated compact metric spaces (which are
not $ANR$'s, still less $CW$ complexes) which are not homotopy equivalent
to a finite $CW$ complex.
\mk 

Wall \cite{Wa2} obtained the finiteness obstruction of a finitely
dominated $CW$ complex $X$ from $C_*(\widetilde X)$, using any 
finite f.g. projective $\bz[\pi_1(X)]$-module chain complex
$${\mathcal P} : \dots \to 0 \to P_n \stackrel{\partial}{\to} P_{n-1} 
\stackrel{\partial}{\to} \dots \stackrel{\partial}{\to}
P_1 \stackrel{\partial}{\to} P_0~.$$
chain equivalent to $C_*(\widetilde{X})$.

\begin{defn} {\normalfont
The {\it projective class} of $X$ is the projective class of $\mathcal P$}
$$[X]~=~\sum\limits^{\infty}_{i=0}(-1)^{i}[P_{i}] \in K_0(\bz[\pi_1(X)])~.
\eqno{\qed}$$
\end{defn}

The projective class is a well-defined chain-homotopy invariant 
of $C_*(\widetilde X)$, with components
$$[X]~=~(\chi(X),[X]) \in K_0(\bz[\pi_1(X)])~=~K_0(\bz)\oplus
\widetilde{K}_0(\bz[\pi_1(X)])~,$$
where 
$$\chi(X)~=~\sum\limits^{\infty}_{i=0}(-1)^i\text{$\#$ of $i$-cells}
\in K_0(\bz)~=~\bz$$
is the Euler characteristic of $X$, and $[X]$ is the
finiteness obstruction.  
\mk

The {\it instant finiteness obstruction} (Ranicki \cite{Ra1})
of a finitely dominated $CW$ complex $X$ is a f.\,g. projective 
$\bz[\pi_1(X)]$-module $P$ representing the finiteness obstruction
$$[X]~=~[P] \in \widetilde{K}_0(\bz[\pi_1(X)])$$ 
which is obtained directly from a 
finite domination $d:K \to X$, $s:X \to K$, a homotopy
$h:d\circ s \simeq \text{id}_X:X \to X$ and the cellular 
$\bz[\pi_1(X)]$-module chain complex $C_*(\widetilde K) $ 
of the cover $\widetilde K=d^*\widetilde X$ of $K$ obtained by pullback from 
the universal cover $\widetilde X$ of $X$, namely
$$P~=~\text{im}(p:\bz[\pi_1(X)]^n \to \bz[\pi_1(X)]^n)$$
with $p=p^2$  the projection 
$$\begin{array}{l}
p~=~\begin{pmatrix}
s \circ d & -\partial & 0 & \dots \\
-s \circ h \circ d & 1- s \circ d & \partial  &\dots \\
s \circ h^2 \circ d & s \circ h \circ d & s \circ d & \dots \\
\vdots & \vdots & \vdots & \ddots
\end{pmatrix}~:\\
\hspace*{20mm} \bz[\pi_1(X)]^n~=~
\sum\limits^{\infty}_{i=0}C_*(\widetilde{K})_i \to 
\sum\limits^{\infty}_{i=0}C_*(\widetilde{K})_i
\end{array}$$
of a f.\,g. free $\bz[\pi_1(X)]$-module of rank 
$$n~=~\sum\limits^{\infty}_{i=0}\text{$\#$ of $i$-cells of $K$}~.$$
In fact, the finiteness obstruction can be obtained in this way
using only the chain homotopy projection 
$q=s\circ d \simeq q^2:C_*(\widetilde K) \to C_*(\widetilde K)$
induced by the homotopy idempotent
$q=s\circ d \simeq q^2 :K \to K$ (L\"uck and Ranicki \cite{LR2}).
\mk 

The finiteness obstruction has many of the usual properties of the
Euler characteristic $\chi$.  For instance, if $X$ is the union of
finitely dominated complexes $ X_1 $ and $ X_2 $ along a common
finitely dominated subcomplex $ X_0 $, then 
\[ [X]=i_{1*}[X_1]+i_{2*}[X_2]-i_{0*}[X_0].\] 
\noindent This is the \textit{sum theorem} for finiteness obstructions,
which was originally proved in Siebenmann's thesis \cite{Si}.

\mk 

The projective class of the product 
$X \times Y$ of finitely dominated $CW$ complexes $X,Y$ is given by
$$[X \times Y]~=~[X] \otimes [Y] \in K_0(\bz[\pi_1(X\times Y)])~,$$
leading to the \textit{product formula} of Gersten \cite{Ge} for the 
finiteness obstruction 
$$[X \times Y] ~=~\chi(X) \otimes [Y] + [X] \otimes
\chi(Y) + [X] \otimes [Y] \in 
\widetilde{K}_0(\bz[\pi_1(X\times Y)])~.$$
In particular, $[X \times S^1]=0$, giving an algebraic proof of
the result (\ref{Mather}) that $X \times S^1$ is homotopy equivalent
to a finite $CW$ complex. 

\mk A fibration $p:E \to B$ with finitely dominated fibre $F$ induces 
{\it transfer maps} in the projective class groups
$$p^{\,!}~:~K_0(\bz[\pi_1(B)]) \to K_0(\bz[\pi_1(E)])~;~[X] \mapsto [Y]$$
sending the projective class of a finitely dominated $CW$ complex $X$
with a $\pi_1$-isomorphism $f:X \to B$ to the projective class of the 
pullback $Y=f^{\,!}E$, which is a finitely dominated $CW$ complex with a
$\pi_1$-isomorphism $f^{\,!}:Y \to E$
\[
\xymatrix@C+4pt{F \ar@{=}[r] \ar[d] &  F\ar[d]  \\
Y \ar[r]^{\displaystyle{f^{\,!}}} \ar[d] & E \ar[d]^{\displaystyle{p}} \\
X \ar[r]^{\displaystyle{f}} & B}
\]
L\"uck \cite{L} obtained the following 
algebraic description of $p^{\,!}$, generalizing
the product formula.\footnote{See
L\"uck and Ranicki \cite{LR1} for the algebraic transfer map in the
surgery obstruction groups.} 
Let $\widetilde{F}$ be the pullback to $F$ of the
universal cover $\widetilde E$ of $E$.
The fibration $p$ induces a morphism of rings
$$U~:~\bz[\pi_1(B)] \to H_0({\rm Hom}_{\bz[\pi_1(E)]}(C_*(\widetilde{F}),C_*(\widetilde{F})))^{op}$$
sending the homotopy class of a loop $\omega:S^1 \to B$ to the chain homotopy
class of the parallel transport chain equivalence 
$U(\omega):C_*(\widetilde{F}) \to C_*(\widetilde{F})$. A f.\,g. projective $\bz[\pi_1(B)]$-module
$$Q~=~{\rm im}(q:\bz[\pi_1(B)]^n \to \bz[\pi_1(B)]^n) \hspace{3mm} (q=q^2)$$
induces a $\bz[\pi_1(E)]$-module chain complex
$$Q^{\,!}~=~{\mathcal C}(U(q):C_*(\widetilde{F})^n \to C_*(\widetilde{F})^n)
\hspace{3mm} (U(q)\simeq U(q)^2)$$
which is algebraically finitely dominated, i.e. chain equivalent to 
a finite f.\,g. projective chain complex. The transfer map is given 
algebraically by
$$p^{\,!}[Q]~=~[Q^{\,!}] \in K_0(\bz[\pi_1(E)])~.$$

\section{The topological space-form problem}

	Another problem in which a finiteness obstruction arises is the 
	\textit{topological space-form problem}.  This is the problem of 
	determining which groups can act freely and properly discontinuously 
	on $S^n $ for some $n$.  
\mk
	
	Swan, \cite{Sw}, solved a homotopy version of this problem by proving 
	that a finite group $ G $ of order $n$ which has periodic 
	cohomology of period $ q $ acts freely on a finite complex of 
	dimension $ dq-1 $ which is homotopy equivalent to a $ (dq-1) 
	$-sphere. Here, $d$ is the greatest common divisor of $n$ 
	and $\phi(n) $, $\phi $ being Euler's $\phi $-function.
	\mk
	One might ask whether such a $ G $ can act on $S^{q-1} $, but 
	this refinement leads to a finiteness obstruction.  It follows from 
	Swan's argument that $ G $ acts freely on a countable $ q-1 
	$-dimensional complex $X$ homotopy equivalent to $S^{q-1} $ 
	and that $ X/G $ is finitely dominated.  The finiteness obstruction 
	of $ X/G $ need not be zero, however, so not every group with cohomology of 
	period $ q $ can act freely on a \textit{finite} complex homotopy 
	equivalent to $S^{q-1} $. Algebraically, the point is that finite 
	groups with $ q $-periodic cohomology have $ q $-periodic 
	resolutions by f.\,g. projective modules but need not 
	have $ q $-periodic resolutions by f.\,g. free modules.
	\mk
	
	After a great deal of work involving both the finiteness obstruction and
	surgery theory, see Madsen, Thomas and Wall \cite{MTW}, it 
	turned out that a group 
	$G$ acts freely on $S^n$ for some $n$ if and only if all 
	of its subgroups of order $p^2$ and $2p$ are cyclic (the condition
	of Milnor).  This 
	is in contrast to the linear case.  A group $ G $ acts linearly on 
	$S^n$ for some $n$ if and only if all subgroups of order 
	$pq$, $p$ and $ q $ not necessarily distinct primes, are 
	cyclic. See Davis and Milgram \cite{DM} for a book-length treatment,
	and Weinberger \cite{We}, p. 110, for a brief discussion.

\section{The Siebenmann end obstruction}

The most significant application of the finiteness obstruction to
the topology of manifolds is via the end obstruction.

\mk

An end $\epsilon$ of an open $n$-dimensional manifold $W$ is {\it tame} if there
exists a sequence $W \supset U_1 \supset U_2 \supset \dots$
of finitely dominated neighbourhoods of $\epsilon$ with
$$\bigcap\limits_i U_i~=~\emptyset~~,~~\pi_1(U_1) \cong \pi_1(U_2)
\cong \dots \cong \pi_1(\epsilon)~.$$
The end is {\it collared} if there exists a neighbourhood of the type
$M \times [0,\infty)$ for some closed $(n-1)$-dimensional manifold $M$,
i.e. if $\epsilon$ is the interior of a compactification $W \cup M$ with
boundary component $M$.

\begin{thm} {\rm (Siebenmann \cite{Si})} A tame end $\epsilon$ 
of an open $n$-dimensional manifold $W$ has an end obstruction
$$[\epsilon] ~=~\underrightarrow{\lim\vspace*{1mm}}_i\,[U_i]
\in \widetilde{K}_0(\bz[\pi_1(\epsilon)])$$
such that $[\epsilon]=0$ if (and for $n\geq 6$ only if) $\epsilon$
can be collared.
\end{thm}

Novikov's 1965 proof of the topological invariance of the rational
Pontrjagin classes made use of the end obstruction in the unobstructed
case when $\pi$ is a free abelian group. The subsequent work of 
Lashof, Rothenberg, Casson, Sullivan, Kirby and Siebenmann
on the Hauptvermutung for high-dimensional manifolds made overt use of the
end obstruction (\cite{Ra5}).

\mk

See sections 7 and 8 below for brief accounts of the applications of
the end obstruction to splitting theorems and triangulation of
high-dimensional manifolds.

\mk 

See Chapman \cite{C2}, Quinn \cite{Q1}, \cite{Q2} and
Connolly and Vajiac \cite{CV} for the controlled end obstruction,
and Ranicki \cite{Ra4} for the bounded end obstruction.

\mk 
See Hughes and Ranicki \cite{HR} for a book-length treatment of ends
and the end obstruction.

\section{Connections with Whitehead torsion}

The finiteness obstruction deals with the existence of a finite $CW$
complex $K$ in a homotopy type, while Whitehead torsion deals with the 
uniqueness of $K$. There are many deep connections between 
the finiteness obstruction and Whitehead torsion, which on the 
purely algebraic level correspond to the connections between the
algebraic $K$-groups $K_0$, $K_1$ (or rather $\widetilde{K}_0$, $Wh$).
\mk

The splitting theorem of Bass, Heller and Swan \cite{BHS}
$$Wh(\pi \times \bz) ~=~Wh(\pi) \oplus \widetilde{K}_0(\bz[\pi]) 
\oplus \widetilde{\text{Nil}}_0(\bz[\pi])\oplus 
\widetilde{\text{Nil}}_0(\bz[\pi])$$
involves a split injection
$$\widetilde{K}_0(\bz[\pi]) \to Wh(\pi \times \bz)~;~
[P] \mapsto \tau(z:P[z,z^{-1}] \to P[z,z^{-1}])~.\eqno{(*)}$$

If $X$ is a finitely dominated space then \ref{Mather} gives a
homotopy equivalence $\phi : X \times S^1 \to K$ to a finite $CW$
complex $K$, uniquely up to simple homotopy equivalence.
Ferry \cite{Fe2} identified the finiteness obstruction $[X] \in
\widetilde{K}_0(\bz[\pi])$ ($\pi=\pi_1(X)$)
with the Whitehead torsion $\tau(f) \in Wh(\pi \times \bz)$ of the
composite self homotopy equivalence of a finite $CW$ complex
$$f~:~K \stackrel{\phi^{-1}}{\longrightarrow} X \times S^1 
\stackrel{\text{id}_X\times -1}{\longrightarrow} X\times S^1
\stackrel{\phi}{\longrightarrow} K~.$$
In Ranicki \cite{Ra2},\cite{Ra3} it was shown that $[X] \mapsto
\tau(f)$ corresponds to the split injection
$$\widetilde{K}_0(\bz[\pi]) \to Wh(\pi \times \bz)~;~
[P] \mapsto \tau(-z:P[z,z^{-1}] \to P[z,z^{-1}])$$
which is different from the original split injection $(*)$
of \cite{BHS}.

\section{The splitting obstruction}

	The finiteness obstruction arises in most classification problems
in high-dimen\-sion\-al topology.  Loosely speaking, proving that two
manifolds are homeomorphic involves decomposing them into homeomorphic
pieces.  The finiteness obstruction is part of the obstruction to
splitting a manifold into pieces.  The nonsimply-connected version of
Browder's $ M \times \br $ Theorem is a case in point.  In \cite{Br},
Browder proved that if $ M^n $, $ n \ge 6$, is a $PL$ manifold without
boundary, $ f:M \to K \times \br^1$ is a ($PL$) proper homotopy
equivalence, and $ K $ is a simply-connected finite complex, then $ M $
is homeomorphic to $ N \times \br^1 $ for some closed manifold $ N $
homotopy equivalent to $ K $.  \mk
	
	\mk
	\centerline{\epsfbox{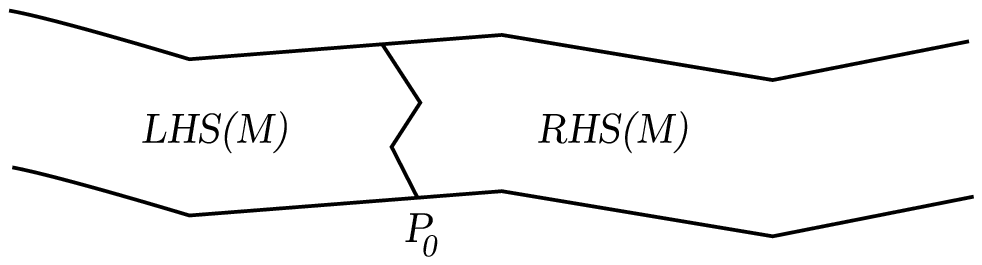}}
	\mk

	
	When $ K $ is connected but not simply-connected, a finiteness 
	obstruction arises.  Here is a quick sketch of the argument: It is not 
	difficult to show that $M$ is 2-ended.  The proper homotopy 
	equivalence $ f:M \to K \times \br^1 $ gives us a proper $PL$ map $ p:M 
	\to \br $.  If $ c \in \br $ is not the image of any vertex, then 
	$ p^{-1}(c) $ is a bicollared $PL$ submanifold of $M$ which 
	separates the ends.  Connected summing components along arcs allows us 
	to assume that $ P_0=p^{-1}(c) $ is connected and a disk-trading 
	argument similar to one in Browder's paper allows us to assume that $ 
	\pi_1P_0 \to \pi_1 M$ is an isomorphism.  See Siebenmann \cite{Si} for 
	details.  An application of the recognition criterion discussed in the 
	third paragraph of this paper shows that the two components of $ M - 
	P_0 $, which we denote by $ RHS(M) $ and $ LHS(M) $, 
	respectively, are finitely dominated.  By the sum theorem, 
	$$ 	[RHS(M)] + [LHS(M)] ~=~0 \in \widetilde{K}_0(\bz[\pi_1(M)])~.$$
	It turns out that the 
	vanishing of $[RHS(M)] = -[LHS(M)] $ is necessary and 
	sufficient for $M$ to be homeomorphic to a product $ N \times \br 
	$, provided that $\dim(M) \ge 6$.  This is one of the main 
	results of \cite{Si}.
	\mk

	It is possible to realize the finiteness obstruction 
	$\sigma\in \widetilde{K}_0(\bz[\pi_1(K)]) $ 
	for an $n$-dimensional manifold $ M^n $ proper homotopy equivalent 
	to $ K \times \br $ for some finite $ K $ whenever $\sigma + 
	(-1)^{n-1}\sigma^{*}=0 $ and $ n \ge 6$.  If we only require that 
	$M$ be properly dominated by some $ K \times \br $, then any 
	finiteness obstruction $\sigma $ can be realized (cf. Pedersen
	and Ranicki \cite{PR}.)	A similar 
	obstruction arises in the problem of determining whether a map $ p: 
	M^{n } \to S^1 $ is homotopic to the projection map of a fiber 
	bundle (Farrell \cite{Fa}).
    
    \mk
     
    The geometric splitting of two-ended open manifolds into right and left sides
    is closely related to the proof of the algebraic splitting theorem
    of Bass, Heller and Swan \cite{BHS} for $Wh(\pi\times \bz)$ -- 
    see Ranicki \cite{Ra4}.

\section{The triangulation of manifolds}

The finiteness obstruction arises in connection with another of 
the fundamental problems of topology: \textit{Is every compact 
topological manifold without boundary homeomorphic to a finite 
polyhedron?} We will examine this problem in much greater detail.
\mk
The triangulation problem was solved affirmatively for two-dimensional 
manifolds by Rado in 1924 and for three-dimensional manifolds by Moise 
in 1952.  Higher dimensions proved less tractable,\footnote{In fact, 
Casson has shown that there are compact four-manifolds without 
boundary which are not homeomorphic to finite polyhedra 
(Akbulut and McCarthy \cite{AM}, p.xvi).
The question is still open in dimensions greater than or equal to five.} a 
circumstance which encouraged the formulation of weaker questions such 
as the following \textit{homotopy triangulation problem}: \textit{Does 
every compact topological manifold have the homotopy type of some 
finite polyhedron?}
\mk
The first solution of this problem came as a corollary to Kirby and 
Siebenmann's theory of $PL $ triangulations of high-dimensional 
topological manifolds.  By a theorem of Hirsch, every topological 
manifold $ M^n $ has a well-defined stable topological normal disk 
bundle.  The total space of this bundle is a closed neighborhood of 
$M$ in some high-dimensional euclidean space.  In \cite{KS1}, Kirby 
and Siebenmann proved that a topological $n$-manifold, $ n \ge 6$, 
has a $PL$ structure if and only if this stable normal bundle 
reduces from $TOP $ to $PL$.  As an immediate corollary, they 
deduced that every compact topological manifold has the homotopy type 
of a finite polyhedron, since each $M$ is homotopy equivalent to 
the total space of the unit disk bundle of its normal disk bundle and 
the total space of the normal disk bundle is a $PL$ manifold because its 
normal bundle is trivial.  The argument of Kirby and Siebenmann also 
shows that each compact topological manifold has a well-defined simple 
homotopy type.  A more refined argument, see p.104 of Kirby and
Siebenmann \cite{KS2}, 
shows that every closed topological manifold of dimension $\ge 6$ 
is a $TOP$ handlebody.  From this it follows immediately that 
every compact topological manifold is homotopy equivalent to a finite 
$CW$ complex and therefore to a finite polyhedron.
\mk
This positive solution to the homotopy-triangulation problem suggests 
that we should look for large naturally-occurring classes of compact 
topological spaces which have the homotopy types of finite polyhedra. 
In 1954, K. Borsuk \cite{Bo2} conjectured that every compact metrizable $ANR$ 
should have the homotopy type of a finite polyhedron. This became 
widely known as \textit{Borsuk's Conjecture}.
\mk
The Borsuk Conjecture was solved by J.  E.  West, \cite{Wes}, 
using results of T.  A.  Chapman, which, in turn, were based on an 
infinite-dimensional version of Kirby-Siebenmann's 
handle-straightening argument.  In a nutshell, Chapman proved that 
every compact manifold modeled on the Hilbert cube ($\equiv 
\prod_{i=1}^{\infty}[0,1]) $ is homotopy equivalent to a finite 
complex and West showed that every compact $ANR$\footnote{A compact 
metrizable space $X$ is an $ANR$ if and only if it embeds as a 
neighborhood retract in separable Hilbert space.  If $X$ has 
finite covering dimension $\le n $, separable Hilbert space can be 
replaced by $\br^{2n+1} $.} is homotopy equivalent to a compact 
manifold modeled on the Hilbert cube. A rather short 
finite-dimensional proof of the topological invariance of Whitehead 
torsion, together with the Borsuk Conjecture was given by Chapman in 
\cite{C1}. See Ranicki and Yamasaki \cite{RY} for a more recent proof,
which makes use of controlled algebraic $K$-theory.
\mk

\providecommand{\bysame}{\leavevmode\hbox to3em{\hrulefill}\thinspace}

\end{document}